\documentclass{amsart}

\usepackage{epsfig,color}

\newcommand{\R}{\mathbb{R}}
\newcommand{\N}{\mathbb{N}}
\newcommand{\C}{\mathbb{C}}
\newcommand{\fc}{1\!\!\;\!{\mbox{l}}}

\theoremstyle{plain}
\newtheorem{lema}{Lemma}[section]
\newtheorem{corolario}[lema]{Corollary}
\newtheorem{proposicao}[lema]{Proposition}

\newtheorem{mainteorema}{Theorem}
\newcommand{\comment}[1]{}

\theoremstyle{definition}
\newtheorem{definicao}[lema]{Definition}
\newtheorem{obs}[lema]{Remark}

\begin{document}
\author{Amanda de Lima and Daniel Smania}

\address{
Departamento de Matem\'atica,
ICMC-USP, Caixa Postal 668,  S\~ao Carlos-SP,
CEP 13560-970
S\~ao Carlos-SP, Brazil}
\email{amandal@icmc.usp.br}
\email{smania@icmc.usp.br}
\urladdr{www.icmc.usp.br/$\sim$smania/}
\date{\today}
\title{On  infinitely cohomologous to zero observables}
\begin{abstract} We show that for a large class of piecewise expanding maps $T$,  the bounded $p$-variation observables $u_0$ that admits an infinite  sequence of bounded $p$-variation observables $u_i$ satisfying
$$u_{i}=  u_{i+1}\circ T-u_{i+1}$$
are  constant.
The method of the proof consists in to find a suitable Hilbert basis for $L^2(hm)$, where $hm$ is the unique absolutely continuous invariant probability of $T$. In terms of this basis,  the action of the Perron-Frobenious and the Koopman operator on $L^2(hm)$ can be easily understood. This result generalizes earlier results  by Bam\'on, Kiwi, Rivera-Letelier and Urz\'ua in the case $T(x)= \ell x \ mod \ 1$, $\ell \in \mathbb{N}\setminus \{0,1\}$ and Lipschitizian observables $u_0$.
\end{abstract}
\thanks{A. de L. is partially supported by FAPESP 04/12117-0 and 2010/17419-6. D.S. is partially supported by CNPq 470957/2006-9, 310964/2006-7 and  303669/2009-8, FAPESP 2003/03107-9,  2008/02841-4 and 2010/08654-1}
\maketitle
\section{Introduction}
Let $T\colon I \rightarrow I$ be a dynamical system. Consider the  {\it cohomological operator} defined by
$$
{\mathcal{L}}(\psi)= \psi \circ T-\psi,
$$
Given an observable, that is, a function $u_0\colon I \rightarrow \mathbb{R}$,  one can
ask if there exists a solution $u_1$ to the {\it Lvisic cohomologous equation}
$${\mathcal{L}}(u_1)=u_0.$$
Such equation was intensively studied after its introduction by
the seminal work of Livsic. These studies  mainly concerns to the
existence and regularity of the solution $u_1$.

Let $\mu$ be an  invariant probability measure  of $T$.
We say that a function $u : I \to \R$ in $L^1(\mu)$ is cohomologous
to zero if there is a function $w:I\to \R$ in $L^1(\mu)$ such that
$$
u = {\mathcal{L}(w)}.
$$
An observable $u_0$ is {\it infinitely cohomologous to zero} if
there exists a  sequence of functions $u_n \in L^1(\mu)$, $n \in
\N$, such that ${\mathcal{L}}^{n}u_n = u_0$, for all $n \in\N$.

Bam\'on, Kiwi, Rivera-Letelier and Urz\'ua \cite{Bamon}
consider the expanding maps defined by
$$
T_\ell (x) = \ell x \ mod \ 1,
$$
where $\ell \ge 2$ is an integer. The  Lebesgue measure on $[0,1]$ is invariant by $T_\ell$. They show that  every
non-constant lipschitzian function $u: I \to \R$  is not
infinitely cohomologous to zero. In this work we generalize this result to a much larger class of observables and piecewise
expanding maps.

In \cite{Bamon} the study of this problem is motivated by the following observation. 
Let  $\lambda \in (-1,1)$, $u_0\colon I \rightarrow \mathbb{R}$ be
a Lipschitz function and define
$$ A\colon I \times \mathbb{R} \rightarrow I \times \mathbb{R}$$
by
$$A_{\lambda, u_0} (x,y)=(T_\ell(x), \lambda y + u_0(x)).$$
In \cite{Bamon} they notice that
\begin{itemize}
\item[i.] If $\mathcal{L}(u_1)=u_0$ then $A_{\lambda, u_1}\circ H=
H\circ A_{\lambda, u_0}$, where $H$ is the homeomorphism
$$H(x,y)=(x,\frac{ y+u_1(x)}{1-\lambda}).$$
\item[ii.] It turns out that the analysis of topological structure
of the attractor of $A_{\lambda,u}$ is  easier if $u$ is {\it not}
cohomologous to zero.\\
\end{itemize}

So if $u_0$ is not infinitely cohomologous to zero, by i. we can reduce the analysis of the topological dynamics of $A_{\lambda,u_0}$ to   the analysis of $A_{\lambda,u_n}$, where $\mathcal{L}^n(u_n)=u_0$ and $u_n$ is not cohomologous to zero.  Using our results, a similar analysis of attractors  could potentially be achieved to  far more general skew-products. 

\subsection{Statement of  results}
Let $I$ be an interval. We say that $T\colon I \rightarrow I$ is a {\it piecewise monotone}  map if there exists a partition by intervals $\{I_1,\dots, I_m\}$ of $I$ such that for each $i\leq m$ the map $T$ is continuous and  strictly  monotone in $I_i$. A piecewise monotone map is {\it onto} if furthermore $T(I_i)=I$ for every $i$. A piecewise monotone map is called {\it expanding} if $T$ is differentiable on each $I_i$ and

$$\inf_i \inf_{x \in I_i} |T'(x)| > 1.$$
In this work, we will consider mainly maps  $T\colon I \rightarrow I$ satisfying  the following conditions:\\
\begin{itemize}
    \item[(D1)] $T$ is piecewise monotone, Lipschitz  on each interval of the partition $I_i$, $i\leq m$. In particular $T'$ is defined almost everywhere and it is an essentially bounded function. We also assume
\begin{equation}\label{inft} ess \inf_m |T'| > 0.\end{equation}
Here $\label{inft} ess \inf_m$ denotes the essential infimum with respect to the Lebesgue measure $m$. 
\item[(D2)] We have  $T(I)=I$ and  moreover  for every interval $H \subset I$ there is a  finite collection of pairwise disjoint open subintervals $H_1,\dots,H_k \subset H$ and $n$ such that $T^{n}$ is a homeomorphism on $H_i$ and
    $$int \  I \subset \cup_i T^n(H_i) .$$
    \item[(D3)] $T$ has a horseshoe, that is, there are  three open intervals  $ J_1, J_2  \subset J \subset I$, with $J_1\cap J_2 = \emptyset$, such that $T$ is a homeomorphism on each  $J_i$ and $T(J_i)=J$, $i=1,2$.
    \item[(D4)] $T$ has an  invariant probability $\mu$ that is absolutely continuous with respect to the Lebesgue measure $m$, so $$\mu(A)= \int_A h \ dm$$
for some $h \in L^1(m)$. We will denote  $\mu=h m$, where $h \in L^1(m)$ and $m$ is the Lebesgue measure on $I$. Moreover  $\mu$ is  exact and  there exist $a,b$ such that
\begin{equation}\label{infh}0< a\leq  h(x) \leq b < \infty \end{equation}
for $hm$-almost every $x$ and the support of $\mu$ is $I$.
\end{itemize}

Our main result is:

\begin{mainteorema}\label{A1}
   Let $T$ be a transformation  satisfying D1-D4 and let  $u_0 \colon I \rightarrow \mathbb{R}$ be an observable with bounded $p$-variation. Then either $u_0$ is constant in $I$ up to a countable set  or  there exist $M \ge 0$ and  bounded $p$-variation functions  $u_i: I \to \R$, with $i\leq M$, which are unique (in $L^1(hm)$ and $BV_{p,I}$) up to an addition by a constant,    such that
 \begin{itemize}
 \item We have
   $${\mathcal{L}}^i u_{i} = u_{0},
   $$
  in $I$ up to a countable set,  for every $i\leq M$.
   \item
  For every  function $\rho$ with bounded
   $p$-variation and every $c \in \mathbb{R}$ we have ${\mathcal{L}}\rho  \neq u_M+c$ in a nonempty open set in $I$. \end{itemize}
\end{mainteorema}

With somehow distinct, but related, assumptions  on $T$ and $u_0$, which are satisfied in  many interesting situations, we can improve this result is such way that ${\mathcal{L}}\rho  \neq u_M+c$ for every $\rho \in L^1(hm)$. In this direction A. Avila \cite{Avila} contributed with improvements of the results in the original version of this work and we are grateful he agreed to include them here. Avila contribution is the following.

\begin{mainteorema}\label{avila2}\cite{Avila} Let $u_0 \in L^1(hm)$ be such that 
$$\int u_0 \ h \ dm=0$$
and such that for every $v \in L^\infty(hm)$  there exist $C> 0$ and $\lambda \in [0,1)$ such that 
$$\big| \int u_0 \cdot  v\circ T^i \cdot  h \ dm  \big| \leq C\lambda^i.$$
Then either $u_0$ is constant $hm$-almost everywhere or  there exist an unique $M \ge 0$ and   functions  $u_i: I \to \R$, with $i\leq M$, $u_i \in L^1(hm)$, which are unique in $L^1(hm)$, up to an addition by a constant,
    such that
    \begin{itemize}
   \item We have $$
   {\mathcal{L}}^i u_{i} = u_{0} \text{ in } L^1(hm)
   $$for every $i\leq M$.
   \item  For every  function $\rho \in L^1(hm)$ and every $c \in \mathbb{R}$ we have ${\mathcal{L}}\rho  \neq u_M+c$ on $L^1(hm)$   \end{itemize}

\end{mainteorema}

Let $(\mathbb{B}, |\cdot  |_\mathbb{B})$ be a Banach space of real-valued, Lebesgue  measurable functions defined on $I$  such that\\
\begin{itemize}
    \item[(D5)] (i) $T$ is a piecewise expanding map satisfying D1 and D4. \\
    (ii)  There exists $C$ and $p_0 \geq 1$ such that $$|f|_{L^1(hm)}\leq C|f|_\mathbb{B}$$ for every $f \in \mathbb{B}$. \\
    (iii) the Perron-Frobenious operator $\Phi_T$ of $T$  is a bounded operator on $\mathbb{B}$ and there exists $h \in \mathbb{B}$, $h >  0$, with $\int h \ dm =1$,  $\lambda \in [0,1)$  and an linear operator $\Psi\colon \mathbb{B} \rightarrow \mathbb{B}$  such that
$$\Phi_T(f) = \int f \ dm \cdot h + \Psi(f),$$
with $$|\Psi^n(f)|_\mathbb{B}\leq C\lambda^n|f|_\mathbb{B},$$
for every $f \in \mathbb{B}$ and $n \in \mathbb{N}$. Moreover  $\Psi(h)=0$. \\
(iv) $1/h \in \mathbb{B}$.\\
(v) The multiplication $$(f,g)\rightarrow f\cdot g$$ is a bounded bilinear transformation on $\mathbb{B}$.\\
(vi) The set $\mathbb{B} $ is dense  in $L^1(hm)$.
\end{itemize}

\begin{mainteorema}\label{A2}
  Let $T$ be a transformation  satisfying D1 and D4 and suppose that  the Banach space of functions $\mathbb{B}$ and $T$ satisfy D5. Let  $u_0 \in \mathbb{B}$ be an observable. Then either $u_0$ is constant $hm$-almost everywhere or  there exist an unique $M \ge 0$ and   functions  $u_i: I \to \R$, with $i\leq M$, $u_i \in L^1(hm)$, which are unique in $L^1(hm)$, up to an addition by a constant,
    such that
    \begin{itemize}
   \item We have $$
   {\mathcal{L}}^i u_{i} = u_{0} \text{ in } L^1(hm)
   $$for every $i\leq M$.
   \item  For every  function $\rho \in L^1(hm)$ and every $c \in \mathbb{R}$ we have ${\mathcal{L}}\rho  \neq u_M+c$ on $L^1(hm)$   \end{itemize}
Moreover $u_i$ belongs to $\mathbb{B}$, for $i\leq M$.
\end{mainteorema}

\begin{obs} In the first version of this work,  Theorem \ref{A2} had additional assumptions. We assumed  for instance that $\mathbb{B}$ was contained in the space of functions with $p$-bounded variation. This is not longer necessary due Avila's contribution (Theorem \ref{avila2}). 

\end{obs}
\begin{obs} The finiteness result for the family of cohomological operators $$\mathcal{L}_\lambda(v) = v\circ T - \lambda v
,$$ with $\lambda \in  (0,1]$, $T(x)= \ell x \mod 1$, for integers $\ell\geq2$ and Lipschitz observables, obtained in \cite[Main Lemma, page 225]{Bamon}, can also  be generalized for  maps described in Remarks \ref{mat}, \ref{bet} and \ref{unimodal}, replacing Lipschitz observables by bounded variation observables. The methods to achieve this generalization are quite similar to those in \cite{Bamon}, so we will not give a full proof here. It is  necessary to use Theorem \ref{A2},  and to replace in their argument the usual Fourier basis by the basis obtained in Section \ref{basis} and the compactness of closed balls centered at zero of the space of Lipschitz functions as  subsets of the space of continuous functions by Helly's Theorem, that is, the compactness of closed balls centered at zero of the space of bounded variation functions as  subsets of $L^1(hm)$.
\end{obs}

\begin{obs}\label{mat} There are plenty of examples of transformations $T\colon I \rightarrow I$ satisfying D1-D4. Let $T$ be a piecewise monotone, expanding map, $C^2$ on each $I_i$. Consider the $m \times  m$ matrix $A_T=(a_{i j})$ defined by $a_{i j}=1$ if
$$\overline{T(int \  I_i)} \subset int \ I_j,$$
and $a_{i j}=0$ otherwise. Here the  closure and interior are taken with respect to the topology of $[0,1]$. Suppose that $A^k_T > 0$ for some $k$. Then $T$  satisfies D1, D2 and D4 and some iteration of $T$ satisfies D1-D4. If we add the assumption that $T$ has a horseshoe, then $T$ fulfills  D1-D4. The space of bounded variation functions $BV(I)$ and $T$ satisfy $D5$.
\end{obs}

\begin{obs}\label{bet} A class of examples satisfying D1-D4 are  $\beta$-transformations $T(x) = \beta x \mod 1$, with $\beta \geq 2$, $\beta \in \mathbb{R}$, $I=[0,1]$. The space of bounded variation functions $BV(I)$ and $T$ satisfy $D5$.
\end{obs}

\begin{obs}\label{unimodal} Let $T\colon [-1,1]\rightarrow [-1,1]$ be a continuous map with $T(-1)=T(1)=-1$, $C^2$ on the intervals $[-1,0]$ and  $[0,1]$, with $T' > 0$ in $[-1,0]$ and  $T' < 0$ in $[0,1]$ and $T(-x)=T(x)$ for every $x \in [-1,1]$. Define
$$\theta = \inf_x |T'(x)|.$$
If $\theta >  1$ then  there exists an unique fixed point $p \in [0,1]$. Define $J=[-p,p]$. If $\theta > \sqrt{2}$ then $T^2$ has a horseshoe in $J$ and satisfies D1-D4 with $I=[T^2(0),T(0)]$. The space of bounded variation functions $BV(I)$ and $T$ satisfies $D5$.
\end{obs}

\begin{obs} Let $T\colon I\rightarrow I$ be a piecewise expanding and  onto map, $C^{1+\alpha_0}$ in each $I_i$, $\alpha_0 \in (0,1)$. Then $T$ satisfies $D1-D4$.  The space of Holder continuous functions $C^\alpha(I)$, for $\alpha\leq \alpha_0$, and $T$ satisfy $D5$.
\end{obs}

\begin{obs} Let $T\colon I\rightarrow I$ be a piecewise expanding  map, linear in each $I_i$. Suppose that $T$ has a horseshoe and satisfies the conditions on the matrix $A_T$ as in Remark \ref{mat}. One can prove using the results of Wong \cite{Wong} that $T$ satisfies $D1-D4$.    The space of bounded $p$-variation functions $BV_{p}(I)$, with $p\geq 1$, and $T$ satisfy $D5$.
\end{obs}

\begin{obs} The mixing assumptions on the  invariant measure $\mu$ are necessary, as it is shown by the following example. Consider a piecewise $C^2$ expanding map $T\colon I \rightarrow I$,  unimodal (continuous and only one turning point), and with a cycle of intervals, that is, there are open intervals $J_j \subset I$, $j< p$ pairwise disjoint, such that $f(\overline{J}_j)\subset \overline{J}_{j+1 \mod p}$ and  $f(\partial J_j)\subset \partial J_{j+1 \mod p}$ . Then $T$ has an absolutely continuous invariant probability $\mu$ and its  support is contained in $\cup_j \overline{J}_j$. Let $\delta \in \mathbb{C}\setminus \{1\}$ be a  $p$-root of unit, $\delta^p=1$.  Define $u_i\colon I \rightarrow \mathbb{C}$, $i\geq 0$, as 
$$u_i(x)= \frac{\delta^j }{(\delta-1)^{i}},  $$
for $x \in J_j$. Define $u_i$ in an arbitrary way elsewhere. It is easy to see that $u_i = u_{i+1}\circ T - u_{i+1}$ on $L^1(hm)$. To obtain real-valued functions, we can consider the real and imaginary parts of $u_i$. 
\end{obs}

\subsection{Topological Results} Replacing Lipschitzian  by bounded p-variation observables has the advantage to allow us to obtain results similar to Theorems \ref{A1} and \ref{A2} to maps which are just {\it topologically conjugate} with maps satisfying the assumptions of those theorems.

We will say that two functions $f,g \colon W\rightarrow \mathbb{R}$ are equal except in a countable subset, $f=g$ on $W$ (e.c.s.) if $\{x \in W \colon \ f(x)\neq g(x)\}$ is countable.

\begin{mainteorema}\label{A1c}
  Let $H\colon I \rightarrow I$ be a homeomorphism, let $T$  be a piecewise monotone map and  $\tilde{T}$ satisfying D1-D4. Suppose that

  $$     H\circ \tilde{T}= T\circ H$$
in  $I$ (e.c.s). Let  $u_0 \colon H(I) \rightarrow \mathbb{R}$ be an observable with bounded $p$-variation. Then either $u_0$ is constant in $H(I)$(e.c.s) or  there exist an unique $M \ge 0$ and bounded $p$-variation  functions  $u_i: H(I) \to \R$,  with $i\leq M$, which are unique up to an addition by a constant (e.c.s.),
    such that
    \begin{itemize}
    \item We have
   $$
   {\mathcal{L}}^i u_{i} = u_{0},
   $$
on $H(I)$ (e.c.s)
   for every $i\leq M$.
   \item  For  every  function $\rho$ with bounded $p$-variation and every $c \in \mathbb{R}$ we have  ${\mathcal{L}}\rho  \neq u_M+c$ in a non-empty open subset in $H(I)$.\end{itemize}
\end{mainteorema}

\begin{mainteorema}\label{A2c}
  Let $H\colon I \rightarrow I$ be a homeomorphism, let $T$  be a piecewise monotone map and  $\tilde{T}$ satisfying D1-D4. Suppose that

  $$     H\circ \tilde{T}= T\circ H$$
in  $I$ (e.c.s.). Suppose that  the  space of functions with bounded $p_0$-variation $BV_{p_0,I}$ and  $\tilde{T}$ satisfy D5. Let  $u_0 \colon H(I) \rightarrow \mathbb{R}$ be an observable with bounded $p_0$-variation. Then either $u_0$ is constant in $H(I)$(e.c.s) or  there exist an unique $M \ge 0$ and continuous (e.c.s.) bounded borelian  functions    $u_i: H(I) \to \R$,  with $i\leq M$, which are unique up to an addition by a constant (e.c.s.),
    such that
    \begin{itemize}
    \item We have
   $$
   {\mathcal{L}}^i u_{i} = u_{0},
   $$
on $H(I)$(e.c.s.)
   for every $i\leq M$.\item  We have  ${\mathcal{L}}\rho  \neq u_M+c$
   \begin{itemize}
   \item[A.] in an  uncountable subset of $H(I)$, if $\rho$  is a  Borel measurable, bounded  function and $c \in \mathbb{R}$.
   \item[B.] in a non-empty open subset of $H(I)$, if $\rho$  is a  Borel measurable, bounded  function  which is continuous in $H(I)$(e.c.s.) and $c \in \mathbb{R}$.
   \end{itemize} \end{itemize}

Moreover $u_i$  has bounded $p_0$-variation, $i\leq M$.
\end{mainteorema}

\begin{obs} Let $T\colon [0,2]\rightarrow [0,2]$ be a piecewise monotone, $C^1$ in $[0,1]$ and $[1,2]$, $T[0,1]=T[1,2]=[0,2]$, with $T(0)=0$, $T'\geq \lambda > 1$ in $[1,2]$ and $T'(x) > 1$ in $x \in (0,1)$ and $T'(0)=1$.  Then $T$ is conjugate with $\tilde{T}(x) = 2\cdot x \mod 1$,  so $T$ satisfies the  assumptions of Theorems \ref{A1c} and \ref{A2c}, considering $p_0=1$ in Theorem \ref{A2c}.
\end{obs}

\begin{obs} Let $T\colon [-1,1]\rightarrow [-1,1]$, $T(-1)=T(1)=-1$, $C^3$ in $[-1,1]$, $T'(0)=0$,  $T' >0$ on $[-1,0)$, $T' < 0$ on $(0,1]$. If $T$ has negative Schwarzian derivative and non-renormalizable then $T$ is conjugate with a tent map
$\tilde{T}_\beta\colon [-1,1]\rightarrow [-1,1]$, defined as $\tilde{T}(x)= -\beta|x|+\beta -1$, with $\beta=exp(h_{top}(T))$. Here  $h_{top}(T)$ denotes the topological entropy of $T$. If $h_{top}(T)\geq \ln(2)/2$ then $T^2\colon I \rightarrow I$, with $I=[T^2(0),T(0)]$, satisfies the assumptions of Theorems \ref{A1c} and  \ref{A2c}, considering $p_0=1$ in Theorem \ref{A2c}.
\end{obs}

\subsection{Continuous observables infinitely cohomologous to zero}  A. Avila told us a nice argument showing the existence of continuous and  non constant observables that are infinitely cohomologous to zero. He kindly agreed to include this result here.

\begin{mainteorema}{\cite{Avila}}\label{cont_ob} Let $T\colon \mathbb{S}^1 \rightarrow \mathbb{S}^1$ be a $C^1$ expanding map on the circle. Then there exists a non constant continuous observable  $u$ that is infinitely cohomologous to zero. 
\end{mainteorema}

\section{Preliminaries}
In this section we present some notations and definitions.

\comment{\begin{definicao}
  We say that $T: I \to I$ is an expanding map if there exists a
  finite partition of the interval $I=[0,1]$, $(I_j)_{j=1, \cdots,
  n}$, such that
  \begin{enumerate}
    \item The restriction of $T$ to each $I_j$, $j = 1, \cdots,
    n$, is of class $C^1$.
    \item The restriction of $T$ to each $I_j$, $j = 1, \cdots,
    n$, is strictly monotonic.
    \item There exists $n_0>0$ a integer such that
    $$
    \inf_{j \in J_{n_0}} \inf_{x_\in {I_j}}|(T^{n_0})'(x)|>1
    $$
    where $(I_j)_{j \in J_{n_0}}$ is the partition of $I$
    associated to $T^{n_0}$.
  \end{enumerate}

\end{definicao}

If $T$ is an $C^2$ expanding and onto map, then $T$ satisfies
$$
\displaystyle{\sup_{j \in J_{n_0}}\sup_{x\neq y \in I_j}
\frac{|(T^{n_0})'(x) - (T^{n_0})'(y)|}{|x-y|}=K< \infty}.
$$

If $T$ is an $C^{1+1/p}$ expanding and onto map,  with $p>1$, then
$T$ satisfies
$$
\displaystyle{\sup_{j \in J_{n_0}}\sup_{x\neq y \in I_j}
\frac{|(T^{n_0})'(x) - (T^{n_0})'(y)|}{|x-y|^{1/p}}=K<\infty}.
$$}

\begin{definicao}
  Given a function $f:I \to \C$  and $p\geq 1$, we define the
  $p$-variation of $f$ by
  $$
  v_{p,I}(f) =\sup \left(\sum_{i=1}^{n} |g(a_i) - g(a_{i-1})|^p\right)^{\frac{1}{p}},
  $$
  where the supremum is taken over all finite sequences   $a_0 < a_1 < \dots < a_n$, $a_i \in I$.

  We say that $f$ has bounded $p$-variation if
  $$v_{p,I}(f)<\infty.
  $$

  Since the Perron-Frobenious operator is not properly defined at points which are image of points where $DT$ is not defined, to define Perron-Frobenious  operator acting in the space of $p$-bounded variation functions it is convenient to identify functions $u$ and $v$ defined on $I$  so that $u=v$  up to a countable subset of $I$. We write $u\sim v$. The set of equivalence classes $[f]$  with respect to the relation $\sim$ such that
  $$v_{p,I}([f])=\inf_{f\sim g} v_{p,I}(g) < \infty$$
  will be called the space of the functions  on $I$ with bounded $p$-variation and denoted
  $BV_{p,I}$. The function $f \rightarrow v_{p,I}([f])$ is a pseudo-norm on $BV_{p,I}$. We can define a norm by
  $$|[f]|_{BV_{p,I}}= \inf_{g\sim f} (\sup |g|  + v_{p,I}(g)).$$
  $(BV_{p,I},|\cdot|_{BV_{p,I}})$ is a Banach space. As usual, from now on we will omit the brackets $[\cdot]$ in the notation of equivalence classes.

  Note that $1/p$-H\"{o}lder continuous functions have
  bounded $p$-variation. When $p=1$, we say that the function has bounded variation.
\end{definicao}

\begin{obs}\label{invhome}
One of the greatest advantages of dealing with $p$-bounded variation observables, in opposition to either H\"{o}lder or Lipschitzian ones, for instance, is that the pseudo-norm $v_{p,I}$ is invariant by homeomorphisms, that is, if $h\colon J \rightarrow I$ is a homeomorphism  and $f\colon I \rightarrow \mathbb{R}$ is an observable  then $$v_{p,I}(f)=v_{p,J}(f\circ h).$$ \end{obs}

\begin{definicao}
  Given a piecewise monotone, expanding map $T$, satisfying $D1$, define the Perron-Frobenius
  operator associated to $T$ by
  $$
  \Phi_Tf(x) = \sum_{j\in J} f(\sigma_j(x))\frac{1}{|T'(\sigma_jx)|}
  \fc_{ T(I_j)} (x),
  $$
  where $\sigma_j: T(I_j)\to I_j$ stands for the inverse branch of
  $T$ restricted to $I_j$ and $\fc_J$ denotes the characteristic function of the set $J$.
\end{definicao}

The main properties of $\Phi_T$ are (see for instance \cite{Broise} and \cite{Baladi}):

$i)$ $\Phi_T$ is a continuous linear operator on $L^1(hm)$.

$ii)$ $\displaystyle{\int_0^1 \Phi_T f \cdot g \, dm = \int_0^1 f
\cdot g\circ T \, dm}$, where $f \in L^1(m)$ and $g \in
L^{\infty}(m)$.

$iii)$ $\Phi_T f = f$ if and only if the measure $\mu = fm$ is
invariant by $T$.

\section{A special basis of  $L^2(hm)$}\label{basis}
In this section we assume that $T$ satisfies $D1$ and $D4$. 
Consider the Hilbert space $L^2(hm)$ with the inner product
$$\left<u, w \right>_{hm}=\int u w h \ dm.$$
Indeed $\left<u, w \right>_{hm}$ is well defined even for $u \in L^k(hm)$ and $w\in L^b(hm)$, with $k, b \in [1,\infty)\cup \{+\infty\}$ satisfying
$$\frac{1}{k} + \frac{1}{b}=1.$$
\noindent Since the measure $h m $ is $T$-invariant we have  
$$\left<u\circ T , w\circ T \right>_{hm}= \left<u, w \right>_{hm}.$$
\noindent In this section we will built a special Hilbert basis for $L^2(hm)$.  Consider the bounded linear operator $P:L^k(hm) \to L^k(hm)$, $k\geq 1$, defined by
$$
P(u) = \frac{\Phi(uh)}{h}.
$$
Due Eq. (\ref{infh}), the operator $P$ is well defined. Indeed 
$$
 \sum_{j\in J} \frac{h(\sigma_j(x))}{h(x)}\frac{1}{|T'(\sigma_jx)|}
  \fc_{ T(I_j)} (x)=1
  $$
  for every $x$ and  $z^k$ is convex, so we have
$$
  \int |Pu|^k h\ dm  \leq \int \big( \sum_{j\in J}\frac{h(\sigma_j(x))}{h(x)} \frac{1}{|T'(\sigma_jx)|} |u|(\sigma_j(x))
  \fc_{ T(I_j)} (x) \big)^k h(x) \ dm
  $$
$$
 \leq \int \sum_{j\in J}\frac{h(\sigma_j(x))}{h(x)} \frac{1}{|T'(\sigma_jx)|} |u|^k(\sigma_j(x))
  \fc_{ T(I_j)} (x)  h(x) \ dm = \int P(|u|^k) h \ dm
  $$
  $$
\leq  \int \Phi(|u|^k)  \ dm =\int |u|^k  \ dm \leq  \frac{1}{a} \int |u|^k  \ hdm.  $$
Note that for $k=1$ we have 
$$ \int |Pu| h\ dm \leq \int \Phi(|u|h)  \ dm=  \int |u|h  \ dm,$$
so $||P||_{L^{1}(hm)}\leq 1$.  

\noindent Let $\mathcal{B}=\{\varphi_i\}_{i \in \N}$ be an orthonormal basis
for  $$\mathrm{Ker}(P)=\{ u \in  L^2(hm) \ s.t. \ P(u)=0   \}.$$ Define
$$
\mathcal{W}=\{ \varphi_i\circ T^j \,\, :\,\, \varphi_i \in
\mathcal{B} \;\; \mbox{e} \;\; j \in \N\}\cup \{\fc_{I}\}.
$$
Recall that $\fc_{A}$ denotes the indicator function of a set $A$. The main result of this section is

\begin{proposicao}\label{lbasel2} Suppose that $T$ satisfies  $D1$ and $D4$. Then 
$\mathcal{W}$ is a Hilbert basis for $L^2(hm)$. Indeed we can choose $\mathcal{B}$ such that $\mathcal{W}\subset L^\infty(hm)$.
\end{proposicao}

\begin{obs}\label{exp_old} A very interesting example of this theorem is given by the function $T\colon [0,1]\rightarrow [0,1]$ defined by  $T(x)= \ell x \mod 1$, with $\ell \in \mathbb{N}\setminus\{0,1\}$. In this case the Ruelle-Perron-Frobenious operator is just
$$(\Phi_T \psi)(x) = \frac{1}{\ell} \sum_{i=0}^{\ell-1} \psi(\frac{x+i}{\ell}).$$
The Lebesgue measure $m$ is an invariant probability, so $P=\Phi_T$. Moreover
$$\mathcal{B}=\{\sin(2\pi n x), \cos(2\pi nx)\colon \ \ell \text{ does not divide $n$}   \}$$
is a basis for $Ker \ P$. Note that
$$\sin(2\pi n T^j(x))= \sin(2\pi n \ell^j x) \ and \  \cos(2\pi n T^j(x))= \cos(2\pi n \ell^j x),$$
so the corresponding set $\mathcal{W}$ is just the classical Fourier basis of $L^2([0,1])$.  \\ \\
\end{obs}

By property ii. of the Perron-Frobenious operator, it is easy to see that the  Koopman operator  $U:L^k(hm) \to
L^k(hm)$, $k\geq 1$, defined by
$$
U(w) = w \circ T,
$$
is the adjoint operator of $P$, that is
\begin{equation} \label{adj}\left<P(u), w \right>_{hm}= \left<u, U(w) \right>_{hm}\end{equation}
for every $u \in L^k(hm)$ and $w\in L^b(hm)$.  Note that $U$ preserves $L^{k}(hm)$ because $hm$ is invariant. Moreover $$P\circ U(f)=f$$ for every $f \in L^1(hm)$.

\begin{lema}\label{ort} $\mathcal{W}$  is an orthonormal set.
\end{lema}
\begin{proof} Indeed
$$|\fc_{I}|_{L^2(hm)}=1,$$
$$|\varphi_i\circ T^j |^2_{L^2(hm)}= |\varphi_i|^2_{L^2(hm)}=1.$$
Futhermore if $$(i_1,j_1)\not = (i_2,j_2)$$
then either $j_1=j_2$, so we have
$$\left<\varphi_{i_1}\circ T^{j_1},\varphi_{i_2}\circ T^{j_2}\right>_{hm}=\left<\varphi_{i_1},\varphi_{i_2}\right>_{hm}=0,$$
or without loss of generality we can assume $j_1 < j_2$ and
$$\left<\varphi_{i_1}\circ T^{j_1},\varphi_{i_2}\circ T^{j_2}\right>_{hm}=\left<\varphi_{i_1},\varphi_{i_2}\circ T^{j_2-j_1}\right>_{hm}=\left<P^{j_2-j_1}(\varphi_{i_1}),\varphi_{i_2}\right>_{hm}=0,$$
and
$$\left<\varphi_{i_1}\circ T^{j_1},\fc_{I}\right>_{hm}= \int  \varphi_{i}\circ T^{j} h \ dm = \int \varphi_{i} h \ dm=\int P( \varphi_{i}) h \ dm =0.$$
\end{proof}

\begin{lema}\label{lema3}
   There exists a countable set
  of functions $\Lambda \subset L^\infty(hm)\cap Ker(P)$  with the following property:  Let $w \in L^k(hm)$, with $k\geq 1$. If  for all $\varphi \in \Lambda$  we have
  $$
  \int w \varphi h \, dm =0,
  $$
  then there exists $\beta \in L^k(hm)$ such that
  $$
  w = \beta \circ T
  $$
$hm$-almost everywhere.
Moreover  $Ker(P)^{\bot}= U(L^2(hm))$.
\end{lema}

\begin{proof}  We claim that  for the existence of $\beta \in L^k(hm)$ such that $w=\beta\circ T$ , it is necessary and sufficient that
  for $hm$-almost every $y \in I$ we have
\begin{equation}\label{eq} \sharp \{ w(x)\colon h(x) \neq 0 \  and \  T(x)=y\} =1.\end{equation}
Indeed,  if the Eq. (\ref{eq}) holds then for every   $y$ satisfying (\ref{eq}),  choosing $x$ such that $T(x)=y$ and $h(x)\neq 0$ we can define
  $$\beta(y)=w(x).$$
If $y$ does not satisfy (\ref{eq}),  define $\beta(y)=0$. Of course $w=\beta\circ T$ $hm$-almost everywhere and, since $hm$ is an invariant measure of $T$, $\beta$   belongs to $L^k(hm)$.

\noindent On the other hand, suppose that  there exists  $\beta \in L^k(hm)$ is such that $w=\beta\circ T$.  Then  $$K= \{ x \colon w(x)=\beta (T(x))\}$$ has full   $hm$-measure. Since  the support of $hm$ is $I$ and $I \subset Im \ T$ it follows that  for $hm$-almost every $y$ we have $\sharp A_y\geq 1$, where 
$$ A_y=\{ w(x)\colon h(x) \neq 0 \  and \  T(x)=y\}.$$
\noindent Suppose there is $\Omega$, with $hm(\Omega) > 0$ such that $\sharp A_y \geq 2$ for every $y \in \Omega$. Note that  $D1$ implies that $f$ and its inverse branches are absolutely continuous functions, so it is easy to see that there are $X_1$ $X_2$ such that $m(X_1), m(X_2) > 0$,  $T(X_i)=\Omega$ and for each $y \in \Omega$ and $i=1,2$ there exists only  one $x_i \in X_i$ such that $T(x_i)=y$. Furthermore $w(x_1)\neq w(x_2)$, $h(x_i)\neq 0$. The absolutely continuity of $T$ and its inverses branches implies that 
$$\tilde{\Omega} = T(X_1\cap K)\cap T(X_2\cap K) \subset \Omega$$
has positive  measure. Let $y\in \tilde{\Omega}$ and $x_i$ as above. Then $w(x_i)=\beta (T(x_i))=\beta(y)$, which contradicts $w(x_1)\neq w(x_2)$.  This concludes the proof of the claim.

\noindent Let $C_i$ be the set of points $x_0 \in I$ such that the function

\begin{equation} \label{diff} F_i(a)=\int_0^a w \circ \sigma_i (Tx)
 \cdot \fc_{T(I_i)}(T(x))\cdot h(x) dm(x)\end{equation}
has derivative $w \circ \sigma_i (T(x_0))\fc_{T(I_i)}(T(x_0))
 h(x_0)$ at $a=x_0$.
 The function in the above integral  belongs to $L^1(m)$, so by the Lebesgue diffentiation theorem the set $$C=\cap_i C_i\setminus \cup_i \partial I_i$$ has full Lebesgue measure in $I$.  Since $T$ is piecewise Lipschitz  we obtain that $$m(T(I \setminus C))=0.$$
Suppose that Eq.  (\ref{eq}) does not hold for $hm$-almost every  $y \in I$. Then   it is not true that Eq.  (\ref{eq})  holds for $hm$-almost every $y \in I \setminus T(I \setminus C).$  Since $hm$-almost every point has at least one preimage $x$ with $h(x)\neq 0$,we conclude that  there exists  $y_0 \in I
\setminus T(I \setminus C)$ and two inverse branches of $T$, denoted by
$\sigma_1$ and $\sigma_2$ such that $y_0$ belongs to the interior of $T(I_1)\cap T(I_2)$ and furthermore
$$
w\circ \sigma_1(y_0) \neq w\circ \sigma_2(y_0), h(\sigma_1(y_0))\neq 0, h(\sigma_2(y_0))\neq 0.
$$
We can assume
$$
w\circ \sigma_1(y_0) > w\circ \sigma_2(y_0),
$$
so
\begin{equation}\label{fdif}
w\circ \sigma_1\circ T\circ \sigma_2(y_0)\fc_{T(I_1)}\circ T\circ \sigma_2(y_0)h\circ \sigma_2(y_0) > w\circ \sigma_2\circ T\circ \sigma_2(y_0)\fc_{T(I_2)}\circ T\circ \sigma_2(y_0)h\circ \sigma_2(y_0) .
\end{equation}

Since  $\sigma_2(y_0) \in C$, the derivatives of the functions $F_1$ and $F_2$ at $a=\sigma_2(y_0)$ are the left  and right hand sides of Eq. (\ref{fdif}) respectively,  so  there exists $\varepsilon>0$ such that
for every closed non degenerate  interval $\tilde{I}_2$ satisfying
\begin{equation}\label{inter} \sigma_2(y_0)\in \tilde{I}_2 \subset (\sigma_2(y_0)-\varepsilon, \sigma_2(y_0)+ \varepsilon) \cap  I_2\end{equation}
we have

$$\int_{\tilde{I}_2} w \circ \sigma_1 (Tx) \fc_{T(I_1)}\circ T(x)\cdot   h(x) dm(x) >  \int_{\tilde{I}_2} w \circ \sigma_2 (Tx) \fc_{T(I_2)}\circ T(x)\cdot  h(x) dm(x).$$

Choose an interval $\tilde{I}_2$ satisfying  Eq. (\ref{inter}) and small enough such that $T(\tilde{I}_2) \subset T(I_1)$. We can assume without loss of generality that $\partial \tilde{I}_2 \subset \mathbb{Q}$.   Then

$$\int_{\tilde{I}_2} w \circ \sigma_1 (Tx) \cdot   h(x) dm(x) >  \int_{\tilde{I}_2} w \circ \sigma_2 (Tx) \cdot   h(x) dm(x).$$

Let  $\tilde{I}_1:= \sigma_1(T(\tilde{I}_2)) \subset I_1$. Define $\varphi$ as

\begin{align}\label{phi}
\varphi(x)= \left\{
\begin{array}{lll}
&- \frac{|T'(x)|}{|T'(\sigma_2(Tx))|}\cdot
\frac{h(\sigma_2(Tx))}{h(x)}&
\;\;\mbox{if} \;\;\; x \in \tilde{I}_1,\\
&1& \;\;\mbox{if} \;\;\; x \in \tilde{I}_2,\\
&0& \;\;\mbox{otherwise.}
\end{array}
\right.
\end{align}

Note that $\varphi \in L^\infty(hm)$ and $\Phi(\varphi h) = 0$.

Hence

\begin{align*}
\int w  \varphi h \, dm &=\\
 & =\int_{\tilde{I}_1} w \varphi h \, dm + \int_{\tilde{I}_2} w \varphi h \, dm\\
&=\int_{\tilde{I}_1}w\cdot \left(-\frac{|T'|}{|T'\circ \sigma_2\circ T|} \frac{h\circ \sigma_2\circ
T}{h}\right) h\, dm + \int_{\tilde{I}_2}w h \, dm.
\end{align*}
Since $\sigma_2 \circ T: \tilde{I}_1 \to \tilde{I}_2$ is Lipschitzian and monotone increasing, we can
make a change of variables to get
\begin{align*}
& - \int_{\tilde{I}_1}w \frac{|T'|}{|T'\circ \sigma_2\circ T|} \frac{h\circ \sigma_2\circ
T}{h} h\, dm
+ \int_{\tilde{I}_2}w  h \, dm \\
& =  - \int_{\tilde{I}_2}w \circ \sigma_1 \circ T \cdot  h \, dm
+ \int_{\tilde{I}_2}w \circ \sigma_2 \circ T \cdot  h \,
dm\\
&< - \int_{\tilde{I}_2}w \circ \sigma_2 \circ T \cdot h \, dm +
\int_{\tilde{I}_2}w \circ \sigma_2 \circ T \cdot h \, dm = 0.
\end{align*}
Therefore
$$
\int w \varphi h \, dm \neq 0.
$$

Let $\Lambda$ be the set of functions $\varphi$ of the form in Eq. (\ref{phi}), with
\begin{itemize}
\item The intervals $\tilde{I}_j \subset I_{i_j}$, $j=1,2$,  and $\sigma_2\colon T(I_{i_2})\rightarrow I_{i_2}$ is the inverse of  $T\colon I_{i_2} \rightarrow T(I_{i_2})$.
\item $T(\tilde{I}_2)=T(\tilde{I}_1)$.
\item $\partial \tilde{I}_2 \subset \mathbb{Q}$.
\end{itemize}

Then it is  easy to see that $\Lambda$ is countable and  $\Lambda \subset L^\infty(hm)\cap Ker P$ and, by the argument above, $\Lambda$ has the wished property.

In particular for $k=2$ we obtain $Ker(P)^{\bot}\subset  U(L^2(hm))$.  The inclusion $U(L^2(hm))\subset Ker(P)^{\bot}$ follows from Eq. (\ref{adj}).
\end{proof}

\begin{proposicao}\label{propfuncao}Let $\Lambda$ be as in Lemma \ref{lema3}.
  Let $u: I \to \R$ be a non constant  function in $L^1(hm)$.
  Then there exists $\varphi \in \Lambda$, and an integer $p \ge 0$
  such that
  $$
  \int u \cdot \varphi \circ T^j \cdot h \, dm = 0, \,\,\, \mbox{for all}
  \,\,\, 0 \le j < p
  $$
  and
  $$
  \int u \varphi \circ T^p \cdot h \, dm \neq 0.
  $$
\end{proposicao}
\begin{proof} Suppose that, for all $\varphi \in \Lambda$ and for all $k \ge 0$
\begin{equation}\label{tk} \int u \, \varphi \circ T^k \cdot h\, dm = 0.\end{equation}
We claim that for every $n$ there exists $\beta_n \in L^1(hm)$ such that
\begin{equation}\label{hypin}
u = \beta_n \circ T^n.
\end{equation}
Indeed, choosing $k = 0$ in Eq. (\ref{tk}) we obtain that  for all $\varphi \in \Lambda$
$$
\int u \, \varphi  h\, dm = 0.
$$
By Lemma \ref{lema3}, there exists $\beta_1 \in L^1(hm)$
such that
$$
u = \beta_1 \circ T.
$$
Suppose by induction that $u = \beta_n \circ T^n$, with $\beta_n \in L^1(hm)$.
By Eq. (\ref{tk}) when $k = n$, for all $\varphi \in \Lambda$ we have
$$\int  \beta_n  \, \varphi  \,h\, dm =\int \beta_n  \circ T^n \cdot  \varphi \circ T^n \cdot h\, dm=\int u \, \varphi\circ T^n \cdot h\, dm=0.
$$
By Lemma \ref{lema3}, there exists $\beta_{n+1} \in L^1(hm)$
such that
$$
\beta_n = \beta_{n+1} \circ T.$$ Hence one has $ u = \beta_{n+1} \circ T^{n+1}.$

Since the measure $hm$ is an exact measure, we can conclude that
$u$ is a constant function. So $u=0$.

\end{proof}

\begin{corolario}\label{cor1}
Let $u:I \to \R$ be a non constant function in $L^2(hm)$. Then there exist
$\varphi_i \in \mathcal{B}$ and an integer $p \ge 0$ such that
$$
\left<u, \varphi_i \circ T^j\right>_{hm} = 0 \;\; \mbox{for all}
\;\; 0\le j <p
$$
and
$$
\left<u, \varphi_i \circ T^p\right>_{hm} \neq 0.
$$
\end{corolario}
\begin{proof} Suppose that for every $\varphi_i \in \mathcal{B}$ and every $j \in \mathbb{N}$
$$
\left<u, \varphi_i \circ T^j\right>_{hm} = 0.
$$
Since $\mathcal{B}$ is a base for  $Ker(P)$ and $U^j\colon L^2(hm)\rightarrow L^2(hm)$ is an isometry,  it follows that
$$\int \varphi \circ T^j \cdot u \cdot h \ dm=0$$
for every $ \varphi \in Ker(P)$ and $j \in \mathbb{N}$. This contradicts Proposition \ref{propfuncao}.
 \end{proof}

\begin{proof}[Proof of Proposition \ref{lbasel2}] It follows from Lemma \ref{ort} and Corollary \ref{cor1} that $\mathcal{W}$  is a basis of $L^2(hm)$.  To construct a basis $\hat{\mathcal{W}} \subset L^\infty(hm)$, consider an enumeration of the set $\Lambda=\{\psi_i\}$ defined in Lemma \ref{lema3}. Apply the Gram–-Schmidt process in the sequence $\psi_i$ to obtain a sequence $\tilde{\psi}_i$ of pairwise orthogonal functions. Discarding the null functions and normalizing the remaining functions, we obtain an orthonormal set of functions $\hat{\mathcal{B}}$. Due Lemma \ref{lema3}
$$\overline{span(\hat{\mathcal{B}})}=\overline{span(\Lambda)}=Ker \ P,$$
so $\hat{\mathcal{B}}$ is a basis of $Ker \ P$, and $$\hat{\mathcal{W}}= \{\phi\circ T^j\colon \ \phi \in \hat{\mathcal{B}}, \ j \in \mathbb{N}   \} \cup \{ \fc_{I} \}$$ is a basis of $L^2(hm)$.   \end{proof}

\begin{corolario}\label{cor_infty}
Let $u:I \to \R$ be a non constant function in $L^1(hm)$. Let $\hat{\mathcal{B}}$ as in the proof of  Proposition \ref{lbasel2}. Then there exist
$\varphi_i \in \hat{\mathcal{B}}$ and an integer $p \ge 0$ such that
$$
\left<u, \varphi_i \circ T^j\right>_{hm} = 0 \;\; \mbox{for all}
\;\; 0\le j <p
$$
and
$$
\left<u, \varphi_i \circ T^p\right>_{hm} \neq 0.
$$
\end{corolario}
\begin{proof} Suppose that for every $\varphi \in \hat{\mathcal{B}}$ and every $j \in \mathbb{N}$
\begin{equation}\label{zero_tudo}
\left<u, \varphi \circ T^j\right>_{hm} = 0.
\end{equation}
Let $\Lambda$ be as in Lemma \ref{lema3}. Since $\hat{\mathcal{B}}$ was obtained applying the Gram–-Schmidt process to $\Lambda$, it follows that Eq. (\ref{zero_tudo}) holds for every $\varphi \in \Lambda$.   This contradicts Proposition \ref{propfuncao}.
 \end{proof}

From now on we assume $\mathcal{W} \subset L^\infty(hm)$. Let  $u \in L^1(hm)$ and consider the Fourier coefficients of $u$ with respect to the basis $\mathcal{W}$

$$
    c_{i,j}(u) =\left< u, U^{j}(\varphi_i)   \right>_{ h m}= \int u \cdot \varphi_i \circ  T^j \cdot  h \, dm.
$$

\begin{proposicao} \label{lemaCoeficiente}
    The  functionals  $c_{i,j}$ have the following properties:
    \begin{enumerate}
        \item $c_{i,j}$ is linear on $L^1(hm)$

        \item $c_{i,j}(U(u)) = c_{i,j-1}(u)$ for $j\geq 1$.

        \item $c_{i,0}(U(u)) = 0$.

        \item $c_{i,j}(P(u)) = c_{i,j+1}(u)$.
    \end{enumerate}
\end{proposicao}

\begin{proof} We have
\begin{enumerate}
    \item The proof is straightforward.

    \item $c_{i,j}(u \circ T) =\left< u \circ T, \varphi_i \circ T^j\right>_{ h m} = \left< u , \varphi_i \circ T^{j-1}\right>_{ h m}= c_{i,j-1}(u).$

    \item  $  c_{i,0}(u \circ T) = \left< U(u), \varphi_i \right>_{  h m} = \left< u, P(\varphi_i) \right>_{  h m} =  \left< u, 0 \right>_{  h m} =0.$
    \item $  c_{i,j}(Pu) = \left< P(u),  U^j(\varphi_i)
        \right>_{h m} = \left< u,  U^{j+1}(\varphi_i)
        \right>_{h m} = c_{i,j+1}(u).$

\end{enumerate}
\end{proof}

\begin{proposicao}\label{l1_lim} For every $u \in L^1(hm)$ and $\varphi_i \in \hat{\mathcal{B}}$ we have
$$\lim_j c_{i,j}=0$$
\end{proposicao} 
\begin{proof} Since $hm$ is exact, it is mixing, so 
$$\lim_j c_{i,j}= \lim_j \int u \cdot \varphi_i \circ  T^j \cdot  h \, dm =0.$$

\end{proof}

\begin{obs} V. Baladi drew to our attention  the method used by M. Pollicott \cite{poli} to built eigenvectors of transfer operators for eigenvalues inside its  essential spectral radius in certain function spaces. In our setting  the method is the following: pick $\varphi \in Ker(P)$ and $|\lambda| < 1$. Then
$$v = \sum_{j=0}^\infty \lambda^j \varphi\circ T^j $$
is a $\lambda$-eigenvector of $P$ in $L^2(hm)$.  Using Propositions \ref{lbasel2} and \ref{lemaCoeficiente} one can easily show that  {\it all } $\lambda$-eigenvalues of $P$ in $L^2(hm)$ , for every $|\lambda|< 1$, can be built in this way. 
\end{obs}

\section{Proof of  Theorem \ref{A1}}

In this section we will study the linear operator
$$
{\mathcal{L}}u = u \circ T - u
$$
acting on functions with bounded $p$-variation  $u: I \to \R$.

First, we will present some properties and then, at the end of
this section, we will prove the theorems announced in
introduction. The following results are well know.

\begin{lema}\label{lema1}
Let ${\mathcal{L}}$ be the linear operator defined above acting on
$L^1(hm)$. Then:
\begin{enumerate}
  \item If $f \in \mathrm{Im}{(\mathcal{L})}$, then $\int f h \, dm =
  0$.
  \item $\mathrm{Ker}({\mathcal{L}}) = \{ f \in L^1(hm)\,\, : \,\, f \,\, \mbox{is constant $hm$-almost
  everywhere}\}$.
\end{enumerate}
\end{lema}

\begin{corolario}\label{uinicidade} Let $u \in L^1(hm)$ and suppose that there exist functions $v, w \in L^1(hm)$ such that
$$\mathcal{L}^n v =u=\mathcal{L}^n w.$$
Then $v = w+ c$ on $L^1(hm)$, for some $c \in \mathbb{R}$. Moreover if  $v, w$ have bounded $p$-variation then
$v = w+ c$ on $I$(e.c.s.).
\end{corolario}
\begin{proof} Define $v_i = \mathcal{L}^{n-i}v$, $w_i = \mathcal{L}^{n-i}w$. We will prove  by induction on $i$ that $v_i=w_i$, if $i < n$ and $v_n = w_n + c$, for some $c \in \mathbb{R}$. Indeed, for $i=0$ we have $w_0=v_0=u$. Suppose that $v_i=w_i$, $i < n$. Then
$$\mathcal{L} (v_{i+1}-w_{i+1})=v_{i}-w_{i}=0,$$
so $v_{i+1}-w_{i+1}$ is $hm$-almost everywhere  constant. If  $i+1 = n$ we are done. If $i+1< n$ then  $\mathcal{L}v_{i+2}=v_{i+1}$ and $\mathcal{L}w_{i+2}=w_{i+1}$, so
$$\int v_{i+1} \ h \ dm = \int w_{i+1} \ h \ dm=0,$$
which implies $c=0$. Now assume that $u,v$ and $w$ have bounded $p$-variation.  Since the support of $hm$ is $I$ and $v=w+c$ $hm$-almost everywhere, we have the $v=w+c$ on a set $\Lambda \subset I$ such that for every non-empty open subset $O$ of $I$ we have that $O \cap \Lambda$ is a dense   and uncountable subset of $O$. Since $v$ and $w$  have just a countable number of discontinuities in $I$, it follows that $v=w+c$ in $I$(e.c.s.).
\end{proof}

\begin{lema}\label{lemacont} Let  $J$ be an {\it open} interval as in D3 and $u \in BV_{p,I}$.
Then  $$v_{p,J}(\mathcal{L}u)\geq v_{p,J}(u)$$ for every $n \in \mathbb{N}$.
\end{lema}
\begin{proof} Let $J_1, J_2 \subset J$ be as in D3. Since $T$ is a homeomorphism on $J_1$ and $J_2$, by  Remark \ref{invhome}
$$v_{p,J}(u\circ T)\geq v_{p,J_1}(u\circ T)+ v_{p,J_2}(u\circ T) = 2 v_{p,J}(u),$$
so
$$v_{p,J}(u\circ T-u)\geq v_{p,J}(u\circ T)-v_{p,J}(u)\geq v_{p,J}(u).$$
\end{proof}

\begin{lema}\label{lema2} There exists $C$ with the following property:
  Let $u_n: I \to \R$, $n\leq M+1$,  be  observables with bounded $p$-variation, $p\geq 1$, such that for every $n\leq M$
  $$
  u_n = {\mathcal{L}}u_{n+1}.
  $$
  Then
  $$
|u_n|_{L^\infty(hm)}\leq  v_{p,I}(u_n) \le C v_{p,I}(u_0)
  $$
 for every $n\leq M$.
\end{lema}
\begin{proof}  Let $J \subset I$ be an one interval as in D3. By Lemma \ref{lemacont}
\begin{equation}\label{boundv} v_{p,J}(u_n)\leq v_{p,J}(u_0)  \end{equation}
for every $n\geq 0$. By D2 there is a  finite collection of pairwise disjoint open intervals $H_1,\dots,H_k \subset J$ and $j$ such that $T^{j}$ is a homeomorphism on each  $H_i$ and
\begin{equation} \label{incl} int \  I \subset \cup_{i=1}^{k} T^j(H_i) .\end{equation}
We claim that for every $\ell \leq j$ and $n$
\begin{equation}\label{indbv} v_{p,T^\ell(H_i)}(u_n)\leq 2^\ell v_{p,J}(u_0)\end{equation}
We will prove this by induction on $\ell$. Of course since $H_i \subset J$, Eq. (\ref{boundv}) implies that for every $i=1,\dots,k$
\begin{equation}\label{boundv2} v_{p,H_i}(u_n)\leq v_{p,J}(u_0),  \end{equation}
So Eq. (\ref{indbv}) holds for $\ell=0$.
Suppose by induction  that Eq. (\ref{indbv}) holds for $\ell < j$ and every $n$. Since $T$ is a homeomorphism on $T^\ell(H_i)$ and $u_{n-1}= u_n\circ T-u_n$ we have
$$v_{p,T^{\ell+1}(H_i)}(u_n)= v_{p,T^\ell(H_i)}(u_n\circ T)\leq v_{p,T^\ell(H_i)}(u_n)+v_{p,T^\ell(H_i)}(u_{n-1})$$
$$\leq 2^{\ell+1} v_{p,J}(u_0).$$
By Eq (\ref{incl})
$$v_{p,I}(u_n)= v_{p,int \ I}(u_n)\leq \sum_{i=1}^k v_{p,T^{j}(H_i)}(u_n)\leq k2^j v_{p,J}(u_0)\leq k2^j v_{p,I}(u_0).$$
Note that since $u_n= u_{n+1}\circ T-u_{n+1}$ it follows that
$$\int u_n h \ dm =0,$$
Suppose that $$ess \  sup_{hm} \ u_n= |u_n|_{L^\infty(hm)}.$$ Then
$$0=\int u_n h \ dm \geq  \ ess \ inf \  u_n = ( ess \ inf \  u_n  - ess \ sup \ u_n )+ ess \ sup \ u_n$$
$$\geq -v_{p,I}(u_n)+ |u_n|_{L^\infty(hm)}.$$
so $|u_n|_{L^\infty(hm)}\leq v_{p,I}(u_n)$. We can obtain the same conclusion for the case $$-ess \  inf_{hm} \ u_n= |u_n|_{L^\infty(hm)},$$ replacing $u_n$ by $-u_n$ in the argument above.
\end{proof}

\begin{proof}[Proof of Theorem \ref{A1}] Define by induction the (either finite or infinite) sequence $u_n\colon I\rightarrow \mathbb{R}$ of functions in the following way: $u_0$ is given. If $u_n$ is defined and there exists a function $v\colon I\rightarrow \mathbb{R}$ with bounded $p$-variation such that $\mathcal{L}v= u_n$ in $L^1(hm)$, then define
$$u_{n+1}= v- \int v \ h \ dm.$$
Otherwise the sequence 	ends with $u_n$.  Note that
$$\mathcal{L}^nu_n=u_0.$$
Define
$$M_0= \sup \{ n \in \mathbb{N}\colon \ u_n \text{ is defined }\} \in \mathbb{N}\cup\{\infty\}.$$
We will show that $M_0< \infty$. Let $M \in \mathbb{N}$, $M\leq M_0$.
Recall the basis $\mathcal{W}$ defined in Section \ref{basis}. By Corollary \ref{cor1} if $u_0$ is not constant almost everywhere there exist $i$ and  $q \ge 0$  such that
$$
c_{i,j}(u_0) = \int u_0 \, \varphi_i\circ T^j\cdot h\, dm = 0,
\;\;\; \mbox{for all} \;\;\; 0 \le j <q
$$
and
$$
c_{i,q}(u_0) = \int u_0 \, \varphi_i\circ T^q\cdot h\, dm \neq 0.
$$
By Lemma \ref{lema2} we have that $|u_n|_{L^{2}(hm)}\leq |u_n|_{L^\infty(hm)} \le C v_{p,I}(u_0)$, so since 
$$|\varphi_i\circ~T^i|_{L^{2}(hm)}~=~1$$ we obtain 
$$
|c_{i,k}(u_n)| = \left|\int u_n \cdot \varphi_i\circ T^k \cdot  h \, dm
\right| \le C 
v_{p,I}(u_0).
$$
Using  Lemma \ref{lemaCoeficiente}, we can now use an argument quite similar to \cite{Bamon}. Observe that
$$
    c_{i,l}(u_{n-1}) = c_{i,l}(u_n \circ T - u_n) = c_{i,l}(u_n \circ T) -
    c_{i,l}(u_n) = c_{i,l-1}(u_n) - c_{i,l}(u_n),
$$
for $l \ge 1$.

For $l=0$,
$$
    c_{i,0}(u_{n-1}) = c_{i,0}(u_n \circ T - u_n) = c_{i,0}(u_n \circ T) - c_{i,0}(u_n)
    = -c_{i,0}(u_n),
$$
for $0 < n \le M$.

Therefore, for $0 < n \le M$
\begin{equation}\label{eqdis1}
    c_{i,l}(u_n) =  c_{i,l-1}(u_n) - c_{i,l}(u_{n-1}), \;
    \mbox{for} \;\; l \ge 1.
\end{equation}

\begin{equation}\label{eqdis2}
c_{i,0}(u_{n-1}) = - c_{i,0}(u_n).
\end{equation}

Since $c_{i,j}(u_0) = 0$ for $0 \le j < q$, by equations
(\ref{eqdis1}) and (\ref{eqdis2}), we can conclude that

\begin{equation}\label{eqdis3}
c_{i,j}(u_n) = 0 \;\; \mbox{for} \;\; 0\le j<q \;\; \mbox{and}
\;\; 0 \le n \le M.
\end{equation}

Now, by equation (\ref{eqdis1}), considering $l = q$, we have
$$
c_{i,q}(u_{n-1}) = c_{i,q-1}(u_n) - c_{i,q}(u_n).
$$

By equation (\ref{eqdis3}), for every $n\leq M$
\begin{equation}\label{eqdis4}
c_{i,q}(u_{n-1}) = -c_{i,q}(u_n).
\end{equation}
By equation (\ref{eqdis4}), we conclude that for $n\leq M$
$$
c_{i,q}(u_n) = (-1)^{n}c_{i,q}(u_0).
$$
Considering  $l = q+1$ in the equation (\ref{eqdis1})
$$ c_{i,q+1}(u_n) =  (-1)^{n}c_{i,q}(u_0) - c_{i,q+1}(u_{n-1})\Rightarrow $$
\begin{equation}\label{tel}  c_{i,q}(u_0) =(-1)^{n}c_{i,q+1}(u_n) + (-1)^{n}c_{i,q+1}(u_{n-1}).\end{equation}
Putting $n=1,\dots,M$ in Eq. (\ref{tel}) and adding the resulting  equations we obtain
\begin{equation}\label{tel2}  M \cdot c_{i,q}(u_0) =(-1)^{M}c_{i,q+1}(u_M) - c_{i,q+1}(u_{0}).\end{equation}
Therefore,
\begin{align*}
M &= \frac{-c_{i,q+1}(u_0) + (-1)^{M}c_{i,q+1}(u_M)}{c_{i,q}(u_0)}  \\
&\le \frac{|c_{i,q+1}(u_0)| + |c_{i,q+1}(u_M)|}{|c_{i,q}(u_0)|}  \\
&\le  \frac{|c_{i,q+1}(u_0)| + Cv_{p,I}(u_0)}{|c_{i,q}(u_0)|}.
\end{align*}
So $M_0$ is  bounded.
Note that by Corollary \ref{uinicidade}, if $v_n \in L^{1}(hm)$ satisfies $\mathcal{L}^nv_n=u_0$ then $v_n=u_n+c$ in $L^1(hm)$, for some $c \in \mathbb{R}$. This proves the uniqueness statements of Theorem \ref{A1}.
\end{proof}
\section{Proof of Theorem \ref{avila2}}

Fix $\lambda < 1$. Denote by $S_\lambda$ the linear  space of the real sequences $x=(x^j)_{j \in \mathbb{N}}$ such that there exists $C $ satisfying
$$|x^j|\leq C \lambda^j.$$
Here we use $x^j$ to denote the $j$-th element of the sequence $x$.  Consider the linear space $\ell_0(\mathbb{N})$ of real sequences  $x=(x^j)_{j\in \mathbb{N}}$ such that 
$$\lim_j x^j=0.$$ 
We define the operator $U\colon \ell_0(\mathbb{N}) \rightarrow \ell_0(\mathbb{N})$ as 
$$U(x)=y,$$
where $y^0 = 0$ and $y^{j+1}=x^j$ for $j \geq 0$. 

We say that $x\in \ell_0(\mathbb{N})$ is infinitely cohomologous to zero with respect to $U$ in $\ell_0(\mathbb{N})$ if there exists an infinite sequence $x_i \in \ell_0(\mathbb{N})$, with $x=x_0$, such that
\begin{equation}\label{coh_abst} x_i = U(x_{i+1})- x_{i+1}.\end{equation}
for every $i\geq0$. 

\begin{lema}\cite{Avila}\label{av2} Let $x \in S_\lambda$. Suppose that there exists a finite sequence $x=x_0$, $x_1, \dots, x_k \in \ell_0(\mathbb{N})$ such that $x_i= U(x_{i+1})-x_{i+1}$ for every $i<k$. Then $x_i \in S_\lambda$, for every $i\leq k$.  If $x$ is infinitely cohomologous to zero with respect to $U$ in $\ell_0(\mathbb{N})$ then $x=0=(0,0,\dots)$. 
\end{lema} 
\begin{proof} Let $x_i \in \ell_0(\mathbb{N})$, $i\leq k$, with $x_0=x$, satisfying Eq. (\ref{coh_abst}) for $i< k$. One can see that 
$$x_{i+1}^j= -\sum_{p\leq j} x_i^p.$$
Since  $\lim_j x_{i+1}^j=0$,  it follows that 
$$\sum_{p} x_i^p =0,$$
consequently since $x_0 \in S_\lambda$ we can prove by induction on $i$  that 
$$|x_{i+1}^j|= |\sum_{p > j} x_i^p| \leq C_i \lambda^j $$
for some $C_i$. We concluded that $x_i \in S_\lambda$ for every $i\leq k$. For each $i\leq k$ we can associate the power series
$$f_i(z)=\sum_{j=0}^\infty x_i^j z^j.$$
Since $x_i=(x_i^j)_j \in S_\lambda$, the power series $f_i$ converges to a complex analytics function on the disc with center at $0$ and radius $1/\lambda > 1$. Note that the sequence $U(x_i)$ is  the sequence of coefficients of the Taylor series (centered at $0$) of the function $zf_i(z)$. So Eq. (\ref{coh_abst}) yields
$$f_i(z)= zf_{i+1}(z)-f_{i+1}(z)= (z-1)f_{i+1}(z)$$
So if $x_0$ is infinitely cohomologous to zero we conclude that
$$f_0(z)=(z-1)^kf_k(z)$$
for every $k$, where $f_k$ is defined in a disc strictly larger than the unit disc.  It follows that $f^{(k)}_0(1)=0$ for every $k$, so $f_0(z)=0$ everywhere. So $x=x_0=0=(0,0,\dots)$. 
\end{proof}

\begin{proof}[Proof of Theorem \ref{avila2}] The Corollary \ref{uinicidade} gives the uniqueness of the sequence $u_i$.   
Now suppose that $u_0$ is infinitely cohomologous to zero. So there exists a sequence $u_i \in L^1(hm)$ such that 
\begin{equation} \label{coh_56}u_i = u_{i+1}\circ T - u_{i+1}.\end{equation}
Consider $\hat{\mathcal{B}}$ as in Corollary  \ref{cor_infty}. Fix $\varphi \in \hat{\mathcal{B}}$. Define the sequence $x_i=(x_i^j)_j$ as 
$$x^j_i = \int u_i \cdot \varphi\circ T^{j} \cdot  h \ dm$$
Since $x_i^j$ are Fourier coefficients of $u_i \in L^1(hm)$ with respect to the Hilbert basis $\mathcal{W}$, by Proposition \ref{l1_lim}  we have
that $\lim_j x^j_i =0$. By Eq. (\ref{coh_56}) and Proposition \ref{lemaCoeficiente} we have 
$$x_i = U(x_{i+1}) - x_i,$$
so $x_0$ is infinitely cohomologous to zero in $\ell_0(\mathbb{N})$. Note that
$$|x_0^j|= \big| \int u_0 \cdot \varphi\circ T^{j} \cdot  h \ dm \big|\leq C \lambda^j, $$
so $x_0 \in S_\lambda$. By Lemma \ref{av2} we have that $x_0=0$. That holds for every $\varphi \in \hat{\mathcal{B}}$, so by Corollary \ref{cor_infty}   the function $u_0$ is zero.

\end{proof}

\section{Proof of Theorem \ref{A2}}\label{imp}

We first make a couple of remarks on condition (D5).  

\begin{obs} Suppose that $T$ and $\mathbb{B}$ satisfy $D5$.  Let $\tilde{h} \in L^1(m)$ be a function satisfying $\Phi_T(\tilde{h})=\tilde{h}$. Then 
\begin{equation}\label{eql1} \tilde{h} = \int \tilde{h} \ dm \cdot   h,\end{equation}
where $h$ is as in D5.iii.   Indeed, by D5.vi there exists a sequence  $h_n \in \mathbb{B}$ such that $h_n\rightarrow_n \tilde{h}$ in $L^{1}(hm)$. Furthermore  since $h$, $1/h \in \mathbb{B}$, due D5.ii there exist $a, b > 0$ such that 
\begin{equation}\label{infhb}0< a\leq  h(x) \leq b < \infty \end{equation}
on $I$. 
So  
\begin{align*} & |\int \tilde{h} \ dm \cdot   h  -  \tilde{h}|_{L^{1}(m)}  \\
&\leq  |\tilde{h}-h_n|_{L^{1}(m)}+ | \int h_n \ dm \cdot   h - \Phi_T^k(h_n) |_{L^{1}(m)}+ | \Phi_T^k(h_n) - \Phi_T^k(\tilde{h}) |_{L^{1}(m)}\\
&\leq 2 |\tilde{h}-h_n|_{L^{1}(m)}+ | \int h_n \ dm \cdot   h - \Phi_T^k(h_n) |_{\mathbb{B}}\\
&\leq 2 |\tilde{h}-h_n|_{L^{1}(m)} +C \lambda^k |h_n|_{\mathbb{B}}. \end{align*}
Given $\epsilon > 0$, choose $n_0$ such that 
$$|\tilde{h}-h_{n_0}|_{L^{1}(m)}\leq \frac{1}{a} |\tilde{h}-h_{n_0}|_{L^{1}(hm)} < \frac{\epsilon}{4},$$
and $k_0$ such that 
$$C \lambda^{k_0} |h_{n_0}|_{\mathbb{B}}< \frac{\epsilon}{2}.$$
Then 
 $$|\int \tilde{h} \ dm \cdot   h  -  \tilde{h}|_{L^{1}(m)} < \epsilon$$
 for every $\epsilon > 0$, so Eq. (\ref{eql1}) holds. In particular if $T$ and $\mathbb{B}$ satisfy D1, D4 and D5 we have those functions $h$ in D4 and D5 coincide. \end{obs}
\begin{obs} \label{relnorm} Note that D5.ii implies that   $$\mathbb{B} \subset  L^1(hm).$$ Moreover D5.iii-v implies that
$$\frac{1}{h}\Psi^j(vh)$$
converges exponentially to zero   in $L^1(hm)$ and  $\mathbb{B}$.
\end{obs}

\begin{lema} \label{rigidity} Let $T$ be a transformation  satisfying D1 and D4 and suppose that  $\mathbb{B}$ and $T$ satisfy D5. Let $u \in \mathbb{B}$ and suppose that there exists $v \in L^1(hm)$ such that
$$u = \mathcal{L} v$$
on $I$. Then $v$ coincides $hm$-almost everywhere with a function $v_1\in \mathbb{B}$.
\end{lema}
\begin{proof} The method we are going to use here is very well known for specific kinds of dynamical systems and observables. See for instance \cite{Broise} for the case of $C^2$ piecewise smooth expanding maps and bounded variation observables.  Replacing $v$ by 
$$v - \int vh \ dm \ \fc_{I},$$  we may assume without loss of generality that $$\int vh \ dm=0.$$ Since
$$u=v \circ T - v,$$
Applying $P^j$, $j\geq 1$, we get
\begin{equation}\label{pe} P^ju=P^{j-1}v  - P^jv,\end{equation}
Putting $j=1,\dots,n$ in Eq. (\ref{pe}) and adding the resulting  equations we obtain
$$v = P^nv + \sum_{j=1}^n P^j u$$
We claim that $|P^jv|_{L^1(hm)}\rightarrow_j 0$. Indeed, due D5.vi for every $\epsilon >0$  there exists $w \in \mathbb{B}$ such that $\int  w \ h \ dm=0$ and $|v-w|_{L^1(hm)}< \epsilon$. Since $||P||_{L^1(hm)}\leq 1$,  for every $j$
$$|P^jv-P^jw|_{L^1(hm)}< \epsilon.$$
Due D5 for every $w \in \mathbb{B}$
$$P^j(w) =  \frac{1}{h}\Psi^j(wh),$$
and  $$|\Psi^j(wh)|_{L^1(hm)} \leq C|\Psi^j(wh)|_\mathbb{B}\leq C\lambda^j|wh|_\mathbb{B},$$ we have that for $j$ large enough
$$|P^jv|_{L^1(hm)} \leq  |P^jv-P^jw|_{L^1(hm)} + |P^jw|_{L^1(hm)}< 2\epsilon.$$
This proves our claim. In particular
$$v = \sum_{j=1}^\infty P^j u,$$
where the convergence of the series is in $L^1(hm)$. On the other hand, by Remark \ref{relnorm} this series converges in $L^1(hm)$ and $\mathbb{B}$ to a function $v_1 \in \mathbb{B}$. So $v=v_1$ $hm$-almost everywhere.
\end{proof}

\begin{proof}[Proof of Theorem \ref{A2}] Since $u_0 \in \mathbb{B}$, by D5,  for every $v \in L^\infty(hm)$ we have
$$\big| \int u_0 \cdot v\circ T^j \cdot  h \ dm \big|=  \big| \int P^j(u_0) \cdot v \cdot  h \ dm \big| \leq C\lambda^j |u_0|_{\mathbb{B}}|v|_{L^\infty(hm)}.$$
By Theorem \ref{avila2} we have that $u_0$ is not infinitely cohomologous to zero in $L^1(hm)$.  Now suppose $\mathcal{L}u_i=u_0$. The uniqueness (up to a constant)  of $u_i$ follows from Corollary \ref{uinicidade}. By Lemma \ref{rigidity} we have $u_i \in \mathbb{B}$. 

  \end{proof}

\section{Topological results}

\begin{proof}[Proof of Theorem \ref{A1c}] Define $\tilde{u}_0=u_0\circ H$. Then $\tilde{u}_0$ has bounded $p$-variation. By Theorem \ref{A1} there exist bounded p-variation functions $\tilde{u}_i$, $i\leq M$, unique up to a constant, such that
$$\tilde{\mathcal{L}}^i \tilde{u_i}=\tilde{u}_0 \ on \ L^1(hm),$$
and  \begin{equation}\label{cont} \tilde{\mathcal{L}}\alpha \neq \tilde{u}_M+c \ on \ L^1(hm),\end{equation} for every bounded $p$-variation function $\alpha$. Here $\tilde{\mathcal{L}}v = v\circ \tilde{T}- v$. Since the support of $hm$ is $I$, it follows that $\tilde{\mathcal{L}}^i \tilde{u_i}=\tilde{u}_0$ in $I$(e.c.s). Define $u_i =\tilde{u}_i\circ H^{-1}$. Then $u_i$ has bounded $p$-variation and
  $$\mathcal{L}^i u_i=u_0 \ on \ H(I) (e.c.s).$$
Suppose that there exists a function $\rho$ with bounded $p$-variation such that $\mathcal{L}\rho= u_M+c$ (e.c.s). Define  $\tilde{\rho}= \rho\circ H$. Then $\tilde{\rho}$ has bounded $p$-variation and $\tilde{\mathcal{L}}\tilde{\rho}= \tilde{u}_M+c$ on $L^1(hm)$. That contradicts Eq. (\ref{cont}). So $\mathcal{L}\rho \neq u_M+c$ in an uncountable subset of $H(I)$. Since the discontinuities of $\mathcal{L}\rho$ and $u_M+c$ are countable, it follows that there is a continuity point $x_0 \in H(I)$ of both functions such that $(\mathcal{L}\rho )(x_0)\neq u_M(x_0)+c$. So there is a non-empty open subset of $H(I)$ such that $\mathcal{L}\rho \neq u_M+c$.
\end{proof}

\begin{proof}[Proof of Theorem \ref{A2c}] The proof of this theorem is quite similar to the proof of Theorem \ref{A1c}. Define $\tilde{u}_0=u_0\circ H$. Then $\tilde{u}_0$ has bounded $p_0$-variation. By Theorem \ref{A2} there exist  bounded $p_0$-variation functions $\tilde{u}_i$, $i\leq M$, unique up to a constant, such that
$$\tilde{\mathcal{L}}^i \tilde{u_i}=\tilde{u}_0 \ on \ L^1(hm),$$
and  \begin{equation}\label{cont2} \tilde{\mathcal{L}}\alpha \neq \tilde{u}_M+c \ on \ L^1(hm),\end{equation} for every  $\alpha \in L^1(hm)$. Here $\tilde{\mathcal{L}}v = v\circ \tilde{T}- v$. Since the support of $hm$ is $I$, it follows that $\tilde{\mathcal{L}}^i \tilde{u_i}=\tilde{u}_0$ in $I$(e.c.s). Define $u_i =\tilde{u}_i\circ H^{-1}$. Then $u_i$ has bounded $p_0$-variation and
  $$\mathcal{L}^i u_i=u_0 \ on \ H(I) (e.c.s).$$

\noindent Now we show  the uniqueness of $u_i$ in the set  of  continuous (e.c.s.), bounded borelian functions.  If continuous (e.c.s.) bounded borelian functions $v_i$ satisfy $\mathcal{L}^iv_i= u_0$ then $\tilde{v}_i= v_i\circ H$ are also continuous (e.c.s.) and moreover   they  belong to $L^1(hm)$ and satisfies $\tilde{\mathcal{L}}^i\tilde{v}_i= \tilde{u}_0$, so by Theorem \ref{A2} we have that  $\tilde{v}_i= \tilde{u}_i +c_i$ for some $c_i \in \mathbb{R}$, where this equality holds in $L^1(hm)$. Since both functions  $\tilde{v}_i$,  $\tilde{u}_i$ are continuous (e.c.s) it follows that  $\tilde{v}_i= \tilde{u}_i +c_i (e.c.s.)$, so  $v_i= u_i +c_i (e.c.s.)$.

\noindent To show conclusions A. and B., suppose that there exists a bounded borelian function $\rho$  such that $\mathcal{L}\rho= u_M+c$ (e.c.s). Define  $\tilde{\rho}= \rho\circ H$. Then $\tilde{\rho}$ is also a bounded borelian function, so it belongs to $L^1(hm)$  and $\tilde{\mathcal{L}}\tilde{\rho}= \tilde{u}_M+c$ (e.c.s), so since $hm$ has  no atoms it follows that this equality holds on  $L^1(hm)$. That contradicts Eq. (\ref{cont2}). So $\mathcal{L}\rho \neq u_M+c$ in an uncountable subset of $H(I)$. If $\rho$ is continuous (e.c.s.) we can now finish the proof exactly as  in the proof of Theorem \ref{A1c}.  \end{proof}

\begin{obs} One can ask why the conclusions of Theorem \ref{A2c} are weaker than those in Theorem \ref{A2}. The problem  is  that the conjugacy between one-dimensional maps  can be singular with respect to the Lebesgue measure. Indeed that is often the case even when the  two one-dimensional maps $T$ and $\tilde{T}$  are very regular, as expanding maps on the circle ( see \cite{ss}). In particular the conjugacy $H$ does not in general preserve either $L^{1}(hm)$, $L^1(m)$ or  the space of Lebesgue measurable functions (see \cite{Goffman}). So note that  if in the proof of Theorem \ref{A2c}  we pick $\rho$ to be either in  $L^1(m)$ or $L^1(hm)$ then it is not true in general that $\rho\circ H$ belongs to $L^1(hm)$. Moreover  since composition with $H$  does not in general preserve Lebesgue measurable functions,  we need to assume that $\rho$ is a Borel measurable function, so $\rho\circ H$ is also Borel measurable. Those are the  reasons why we assume that $\rho$ is bounded  and borelian in Theorem \ref{A2c}. 
\end{obs}

\section{Observables infinitely cohomologous to zero}

Consider the Banach space of summable sequences $\ell^1(\mathbb{N})$. For a sequence $x=(x^j)_{j\in \mathbb{N}}$ denote 
$$|x|_{\ell^1(\mathbb{N})}= \sum_j |x^j|.$$
We define the operator $U\colon \ell^1(\mathbb{N}) \rightarrow \ell^1(\mathbb{N})$ as the norm preserving map
$$U(x)=y,$$
where $y^0 = 0$ and $y^{j+1}=x^j$ for $j \geq 0$.

We say that $x\in \ell^1(\mathbb{N})$ is infinitely cohomologous to zero with respect to $U$ if there exists an infinite sequence $x_i \in \ell^1(\mathbb{N})$, with $x=x_0$, such that
$$x_i = U(x_{i+1})- x_{i+1}.$$
for every $i\geq0$.

\begin{lema}{\cite{Avila}} \label{avila} There is a non vanishing sequence $x \in \ell^1(\mathbb{N})$ which is infinitely cohomologous to zero with respect to $U$. 
\end{lema}
\begin{proof} We claim that for every $k \in \mathbb{N}$ there exist $$x_{0,k}, \ x_{1,k}, \ \dots,\  x_{k,k} \in \ell^1(\mathbb{N}),$$ all of them with compact support,  such that $x_{0,k}^0=1$, $$x_{i,k} = U(x_{{i+1},k})-x_{i+1,k}$$
\begin{equation} \label{cauc} |x_{i,k+1}-x_{i,k}|_{\ell^1(\mathbb{N})} < 2^{-k-1},\end{equation}
for every $i< k$.  

The proof is by induction on $k$. Choose $x_{0,0}=(1,0,0,0,\dots)$.  Suppose by induction we found a finite sequence  $x_{i,k}$, $i\leq k$, with the properties above. Fix $N > 0$. Define $x_{k,k+1}$ as $x_{k,k+1}^0=x_{k,k}^0$,
$x_{k,k+1}^j= x_{k,k}^j-\delta/N$, for $1\leq j\leq N$, and $x_{k,k+1}^j= x_{k,k}^j$ for $j\geq N+1$. Here $\delta=\sum_j x_{k,k}^j$.  Defining $$x_{k+1,k+1}^j = -\sum_{p \leq j} x_{k,k+1}^p,$$
we have that $x_{k+1,k+1}$ has compact support and $x_{k,k+1}= U(x_{k+1,k+1})- x_{k+1,k+1}$. Now define by induction
$$x_{i,k+1}= U(x_{i+1,k+1}) - x_{i+1,k+1}, \ i < k.$$
In particular $x_{i,k+1}^0= -x_{i+1,k+1}^0$ for $i\leq k$. Since  $x_{i,k}^0= -x_{i+1,k}^0$ for $i<  k$ and  $x_{k,k+1}^0=x_{k,k}^0$ we have $x_{0,k+1}=1$.  Furthermore it is not difficult to see that if $N$ is large enough then $$|x_{i,k+1}-x_{i,k}|_{\ell^1(\mathbb{N})} < 2^{-k-1},$$
for every $i< k$. This completes the inductive step.

By Eq. (\ref{cauc}), for every $i$  there exists  $x_{i} \in \ell^1(\mathbb{N})$ such that $\lim_{k} x_{i,k} = x_{i}$ on $\ell^1(\mathbb{N})$.  It is easy  to check that $x_i= U(x_{i+1})- x_{i+1}$ and $x_0^0=1$. Pick $x=x_0$. \end{proof}

\begin{proof}[Proof of Theorem \ref{cont_ob}] Since $T$ is topologically conjugate with $T_\ell = \ell x \mod 1$, $\ell \in \mathbb{Z}\setminus \{-1,0,1\}$, it is enough to show the Theorem  \ref{cont_ob} for $T_\ell$.   Choose $n$ such that $\ell$ does not divide $n$. Let $x=(x_j)_j \in \ell^1(\mathbb{N})$ as in Lemma \ref{avila}. Define
$$u(x)= \sum_{j=0}^\infty x_j\sin(2\pi n\ell^j \ x).$$
The function $u$ is continuous and non constant. Using Remark \ref{exp_old} and Proposition \ref{lemaCoeficiente} one can easily show that $u$ is infinitely cohomologous to  zero. 
\end{proof} 

\section*{Acknowledgment} We are especially grateful to A. Avila for his contributions to  this work. We also would like to thank   V. Baladi, A. Lopes, J. Rivera-Letelier,  A. Tahzibi, A. Wilkinson and the referee for the very useful comments and suggestions.

 \end{document}